\documentstyle[11pt,amssymb]{article}
\newcounter{SNO}
\newcounter{THNO}[section]
\renewcommand{\theSNO}{A.\arabic{SNO}}
\renewcommand{\theTHNO}{\arabic{section}.\arabic{THNO} }

\def\sec#1{\refstepcounter{section}\par\vspace{0.9cm}
\par\noindent
{\large\bf\arabic{section} \ #1    }\nopagebreak\par\vspace{0.3cm}}
\def\sect#1{\par\vspace{0.9cm}
\par\noindent
{\large\bf \ #1    }\nopagebreak\par\vspace{0.3cm}}

\def\ap#1{\refstepcounter{SNO}\par\vspace{0.5cm}
\par\noindent\begingroup \it
\leftskip=0em\hspace{0em}{\bf\ #1 \theSNO\ }}
\def\eap{\par\endgroup}

\def\th#1{\refstepcounter{THNO}\par\vspace{0.5cm}
\par\noindent\begingroup \it
\leftskip=0em\hspace{0em}{\bf\ #1 \theTHNO\ }}
\def\eth{\par\endgroup}
\def\pr{\par\noindent{\bf\ Proof. }}
\def\db#1{ D^b_{coh}({#1})}
\def\d#1{ {\cal D}^b ({#1})}
\def\h#1,#2{{\rm Hom}({#1}\:,\; {#2})}
\def\H#1,#2,#3,#4{{\rm Hom}^{#1}_{#2}({#3}\:,\; {#4})}
\def\E#1,#2,#3,#4{{\rm Ext}^{#1}_{#2}({#3}\:,\; {#4})}

\def\lto{\longrightarrow}
\def\o#1{{\cal O}_{#1}}

\def\ss#1{\scriptsize{$#1$}}
\def\da{\big\downarrow}
\def\ua{\big\uparrow}
\def\dda{\begin{picture}(4,15)\multiput(2,9)(0,-2){5}{\circle*{1}}
\put(2,1){\vector(0,-1){4}}\end{picture}}

\title{\bf Reconstruction of a variety from the derived category and  groups of autoequivalences.}
\author{A.~Bondal \and D.~Orlov}
\date{}
\begin{document}

\maketitle
%
%
%
%
%
%
%
{\large \bf Introduction.}
\vspace{0.3cm}

There exist examples of different varieties $X$ having
equivalent the derived categories ${\db{X}}$ of coherent sheaves.
For abelian varieties and K3 surfaces this kind of equivalences were constructed by Mukai, Polishchuk and the second author in \cite{Mu1}, \cite{Mu2}, \cite{Pol}, \cite{Or}. In \cite{BO} we prove 
equivalence of the derived categories for varieties connected by some kind of flops.

 Does it mean that
$\db{X}$ is a weak invariant of a variety $X$? In this paper we  shall show that this is not the case at least for some types of algebraic varieties.

We prove that a variety $X$ is uniquely determined by its category ${\db{X}}$,
 if its anticanonical (Fano case) or canonical (general type case) sheaf
is ample. 

To reconstruct a variety from the category we, in fact, use nothing but a  graded structure of the category, i.e. we only need to fix the translation functor.

 The idea is
that for good, in the above sense, varieties we can recognize the skyscraper sheaves of closed points in $\db{X}$. The main tool for this is the Serre
functor \cite{BK} (see also ch.1), which for ${\db{X}}$ can be regarded as a categorical incarnation of the canonical sheaf $\omega_X$.

With respect to the above problem it is natural to introduce the groupoid with object being the categories ${\db{X}}$ and with morphisms being autoequivalences.

There are two natural questions related to a groupoid: which objects are isomorphic and what is the group of automorphisms of an individual object. The first problem was considered above in the framework of graded categories. To tackle the second one we need triangulated structure of the category. In chapter 3 we prove that for a smooth algebraic variety with ample either canonical or anticanonical sheaf the group of exact autoequivalences is the semidirect product of the group of automorphisms of the variety and the Picard group plus translations.

The answers to the above questions for the case of varieties with non-ample and non-antiample canonical sheaf seem to be of considerable
interest.

We are thankful to H.Meltzer for a stimulating discussion and to Max-Plank-Institut f\"ur Mathematik for hospitality.
\sec{Preliminaries.}

We collect here some facts relating to functors in graded and triangulated
categories, especially those about Serre functor.

In this paper for simplicity reasons we consider only $k$-linear additive
categories, where $k$ is an arbitrary field.

By definition a {\sf graded category} is a couple $( {\cal D}, T_{\cal D} )$ consisting of a category
${\cal D}$ and a fixed equivalence functor $T_D : {\cal D}\lto {\cal D}$ , called translation
functor.

Recall that a {\sf triangulated category} is a graded category with
an additional structure:  a  distinguished class of exact triangles
satisfying certain axioms (see \cite{ver}).

A functor $F : {\cal D} \lto{\cal D}'$ between two graded categories
${\cal D}$ and ${\cal D}'$ is called {\sf graded} if it commutes with the translation functor. More precisely, there is fixed a natural isomorphism of functors:
$$
t_F : F\circ T_{\cal D} \stackrel{\sim}{\to}T_{\cal D'} \circ F.
$$

In the sequel we omit subscripts in the notation of translation functors
because from their position in formulas it is always clear to which category
they belong.

While considering graded functors, we use graded natural transformations.
A natural transformation $\mu$ between graded functors $F$ and $G$ is called {\sf graded} if the following diagram is commutative:
$$
\begin{array}{ccc}
F\circ T& \stackrel{t_F}{\to}& T\circ F\\
\da\rlap{\ss{\mu T}}&&\da\rlap{\ss{T\mu}}\\
G\circ T& \stackrel{t_G}{\to}&T\circ G.
\end{array}
$$

A graded functor $F : {\cal D}\lto{\cal D}'$ between triangulated categories
is called {\sf exact} if it  transforms all exact triangles into exact
triangles in the following sense. If $X\to Y\to Z\to TX$ is an exact triangle
in ${\cal D}$, then one takes $FX\to FY\to FZ\to FTX$ and substitutes in this
 sequence $FT(X)$ by $TF(X)$ using the natural isomorphism of $FT$ with $TF$.
The result
$$
FX\to FY\to FZ\to TFX
$$
should be an exact triangle in ${\cal D}'$.

A morphism between exact functors is, by definition, a graded natural
transformation.

A functor, which is isomorphic to an exact functor, can be endowed with
a structure of a graded functor so that it becomes an exact functor. Indeed,
 if $F$ is exact , then using isomorphism $\mu : F\stackrel{\sim}{\to} G$
one constructs the natural isomorphism $t_G : GT\stackrel{\sim}{\to}TG$,
$t_G=\mu t_F\mu^{-1}$, which makes $G$ graded.
Since any triangle , isomorphic to an exact triangle, is again exact, then $G$ transforms exact triangles into exact ones.
The natural transformation $\mu$ becomes a graded transformation of exact functors.

Let $F : {\cal D}\lto{\cal D}'$ be a functor. Suppose we fix a class ${\cal C}$ of objects in ${\cal D}$ and for any object $X\in {\cal C}$ some object $X'$ isomorphic to $FX$ in ${\cal D}'$. If we additionally fix for any $X\in {\cal C}$ an isomorphism $FX\stackrel{\sim}{\to}X'$ then there exists a new functor $G : {\cal D}\to {\cal D}'$, which is isomorphic to $F$ and such that

\begin{eqnarray}\label{cor}
GX=&FX,& \mbox{for}\; X\notin{\cal C},\nonumber\\
GX=&X',& \mbox{for}\; X\in{\cal C},
\end{eqnarray}
with the evident action on morphisms.

We shall frequently use this simple fact in the sequel.

\th{Proposition}\label{adjf}
i) Let $F : {\cal D} \lto{\cal D}'$ be  a graded functor between graded
categories, $G : {\cal D}' \lto{\cal D}$ its left adjoint, so that the natural
transformations are given:
\begin{equation}\label{isom}
id_{\cal D'}\stackrel{\alpha}{\lto}F\circ G , \quad G\circ F\stackrel{\beta}{\lto} id_{\cal D}.
\end{equation}

Then $G$ can be canonically endowed with the structure of a graded functor,
such that (\ref{isom}) become morphisms of graded functors,

ii) if, in addition, $F$ is an exact functor between triangulated categories,
then $G$ also becomes an exact functor.
\eth
\pr
Let us make $G$ graded.

By the adjointness of $G$ and $F$ and since $T_{\cal D}$ and $T_{\cal D'}$
are equivalences we have the following sequence of bifunctorial isomorphisms:
\begin{eqnarray}\label{adj}
&{\h {GTX}, Y}\cong{\h TX, {FY}}\cong {\h X, {T^{-1} FY}}\cong&\nonumber\\
&{\h X, {FT^{-1}Y}}\cong {\h GX, {T^{-1}Y}}\cong{\h TGX, Y}&
\end{eqnarray}
for any $X\in {\cal D}', Y\in{\cal D}$.

By the well known Brown lemma \cite{Br} this gives a functorial isomorphism:
$$
t_G : GT\lto TG.
$$

Taking $Y=TGX$ in (\ref{adj}) and tracking carefully
the preimage in ${\h GTX, {TGX}}$ of $id_{TGX}$ in ${\h TGX, {TGX}}$ under
the chain of isomorphisms in (\ref{adj}) one obtains a formula for $t_G$.
It is, in fact, canonically given as the composition of the following
sequence of natural transformations:
\begin{equation}\label{GTG}
GT\stackrel{GT\alpha}{\lto}GTFG\stackrel{Gt^{-1}_F G}{\lto}
GFTG\stackrel{\beta TG}{\lto}TG.
\end{equation}

Here we use morphisms $\alpha$ and $\beta$ from (\ref{isom}) and the grading
isomorphism $t_F$ for $F$:
$$
t_F : FT\lto TF.
$$

To show that, say, $\alpha$ is an isomorphism of graded functors is equivalent
to prove that the diagram
$$
\begin{array}{ccccc}
T&\stackrel{T\alpha}{\lto}&TFG&\stackrel{t_F G}{\lto}& FTG\\
\llap{\ss{\alpha T}}\da&&\llap{\ss{\alpha TFG}}\dda&&\llap{\ss{\alpha FTG}}\dda
\ua\rlap{\ss{F\beta TG}}\\
FGT&\stackrel{FGT\alpha}{\lto}&FGTFG&\stackrel{FGt_F G}{\lto}&FGFTG,
\end{array}
$$
being considered without dotted arrows, is commutative.
One can split it by dotted arrows into two commutative squares and the loop,
the latter being commutative due to the fact that for adjoint functors the
composition
$$
F\stackrel{\alpha F}{\lto}FGF\stackrel{F\beta}{\lto}F
$$
equals $id_F$.

Notice that the inverse morphism to (\ref{GTG}) is given by the composition
\begin{equation}\label{TGT}
TG\stackrel{TGT^{-1}\alpha T}{\lto}TGT^{-1}FGT\stackrel{TGT^{-1}t_F T^{-1}GT}
{\lto}TGFT^{-1}GT\stackrel{T\beta T^{-1}GT}{\lto}GT.
\end{equation}

That can be found in the same way as (\ref{GTG}) by putting $Y=GTX$ in (\ref{adj}).

One needs a lot of commutative diagrams to prove directly, without use of
(\ref{adj}), that (\ref{GTG}) and (\ref{TGT}) are mutually inverse.

ii) \cite{BK} Let  $A\stackrel{\alpha}{\to}B\to C\to TA$ be an exact triangle
 in ${\cal D}'$. We have to show that $G$ transforms this exact triangle
into an exact one.

Let us insert the morphism $G(\alpha) : GA\to GB$ into an exact triangle:
$$
GA\to GB\to Z\to TGA
$$
Applying functor $F$ to it we obtain an exact triangle:
$$
FGA\to FGB\to FZ\to TFGA
$$
(we use henceforth with no mention the commutation isomorphisms like
$TF\stackrel{\sim}{\to}FT$).

By means of $id{\to}FG$ we have a commutative diagram:
$$
\begin{array}{ccccccc}
A & \lto & B & \lto & C & \lto & TA\\
\da&&\da&&&&\da\\
FGA & \lto & FGB & \lto & FZ & \lto & TFGA
\end{array}
$$

By axioms of triangulated categories there exists a morphism $\mu :C\to FZ$
that completes this commutative diagram. By adjunction we obtain a morphism
$\nu : GC\to Z$
that makes commutative the following diagram:
$$
\begin{array}{ccccccc}
GA & \lto & GB & \lto & GC & \lto & TGA\\
\llap{\ss{id}}\da &&\llap{\ss{id}} \da&&\llap{\ss{\wr}}\da\rlap{\ss{\nu}} &&\llap{\ss{id}}\da\\
GA & \lto & GB & \lto & Z & \lto & TGA
\end{array}
$$
Therefore the lower triangle is exact.

If $F$ is a graded autoequivalence in a graded category, then the adjoint
functor is its quasi-inverse.

We may consider a category with object being graded (or respectively triangulated) categories and morphisms being isomorphic classes of graded (respectively exact) equivalences. The proposition ensures that this category is a groupoid. In particular the set of isomorphism classes
of graded autoequivalences in a graded category or of exact autoequivalences
in a triangulated category is a group.

Now we outline main properties of the  Serre functor. Its abstract definition was introduced in \cite{BK}.
\th{Definition}\label{SF}
Let ${\cal D}$ be a $k$-linear category with finite--dimensional
${\rm Hom's}$. A covariant additive  functor $S: {\cal D}\to{\cal D}$
is called a {\sf Serre functor} if it is a category equivalence and there
are given bi--functorial isomorphisms
$$
\varphi_{A,B}: {\rm Hom}_{\cal D}(A\:, \;B)\stackrel{\sim}{\lto}{\rm Hom}_{\cal D}(B\:, \;SA)^*
$$
for any $A,B\in{\cal D}$, with the property that the following diagram is commutative:
$$
\begin{array}{ccc}
 {\rm Hom}_{\cal D}(A\:, \;B) & \stackrel{\varphi_{A,B}}{\lto} & {\rm Hom}_{\cal D}(B\:, \;SA)^*\\
\llap{\ss{\wr}}\da && \llap{\ss{\wr}}\ua\\
{\rm Hom}_{\cal D}(SA\:, \;SB) & \stackrel{\varphi_{SA,SB}}{\lto} & {\rm Hom}_{\cal D}(SB\:, \;S^{2}A)^{*}
\end{array}
$$
The vertical isomorphisms in this diagram are induced by $S$.
\eth

\th{Proposition}\label{com}
Any autoequivalence $\Phi :{\cal D}\lto{\cal D}$ commutes with a Serre
functor, i.e. there exists a natural graded isomorphism of functors
$\Phi\circ S\stackrel{\sim}{\to}S\circ\Phi$.
\eth
\pr
For any couple of objects $A, B$ in ${\cal D}$ we have a system of natural
isomorphisms:

\begin{eqnarray}\label{cong}
&&{\h \Phi A,{ \Phi SB}}\cong{\h A,{ SB}}\cong{\h B, A}^{*}\cong{\h \Phi B, {\Phi A}}^{*}\cong\nonumber\\
&&\cong{\h \Phi A, { S\Phi B}}
\end{eqnarray}
Since $\Phi$ is an equivalence the essential image of $\Phi$ covers the whole
${\cal D}$, i.e. up to isomorphism any object can be presented as $\Phi A$
for some $A$. This means that (\ref{cong}) gives isomorphism of the
contravariant functors represented by objects $ \Phi SB$ and $S\Phi B$. By
the Brown lemma \cite{Br} morphisms  between representable functors are in
 one-one correspondence with those between the representation objects. This
gives isomorphism
$$
\Phi SB\stackrel{\sim}{\to}S\Phi B,
$$
which is, in fact, natural with respect to $B$.

\th{Proposition}\label{GE}
i) Any  Serre functor in a graded category is graded;

ii) A Serre functor in a triangulated category is exact.
\eth
\pr
i) follows from previous proposition.

ii) the fact that a Serre functor takes exact triangles into exact ones
is proved in \cite{BK}.

\th{Proposition}\label{USF}{\rm \cite{BK}}
Any two Serre functors are connected by a canonical graded
functorial isomorphism, which commutes with the bifunctorial isomorphisms
$\phi_{A, B}$ in the definition of Serre functor.
\eth
\pr
Let $S$ and $S'$ be two Serre functors in a  category ${\cal D}$. Then for
any object $A$ in ${\cal D}$ we have natural isomorphisms:
$$
{\h A, A}\cong{\h A, {SA}}^*\cong{\h SA, {S'A}}
$$
Taking the image of the identity morphism $id_A$ with respect to this identification
 we obtain a morphism $SA\to S'A$, which, in fact, gives a graded functorial
isomorphism $S\stackrel{\sim}{\to}S'$, which commutes with $\phi_{A, B}$.

Thus, a Serre functor in a category ${\cal D}$, if it exists, is unique up
to a graded natural isomorphism. By definition it is intrinsically related
to the structure of the category. We shall use this later to reconstruct a
variety from its derived category and to find the group of exact
autoequivalences for algebraic varieties with ample either canonical or
anticanonical sheaf.

\sec{Reconstruction of a variety from the derived category of coherent sheaves.}

In this chapter we show that a variety $X$ can be uniquely reconstructed from
the derived category of coherent sheaves on it, provided $X$ is smooth and has
ample either canonical or anticanonical sheaf. From the category we need only
its grading, i.e. fixation of the translation functor.

Roughly, the reconstruction  proceeds as follows.
First, by means of the Serre functor we distinguish the skyscraper sheaves
of closed points in the variety. Then we find the invertible sheaves and use them to define the Zarisky
topology and the structure sheaf  for the variety.

Let ${\cal D}$ be a k--linear category. Denote by $S_D$ the
Serre functor in ${\cal D}$ (for the case it exists).

Let $X$ be a smooth algebraic variety, $n={\rm dim}X$, ${\cal D}=\db{X}$
the derived category of coherent sheaves on $X$ and $\omega_X$ the
canonical sheaf. Then the functor
\begin{equation}\label{Ser}
(\cdot)\otimes\omega_X[n]
\end{equation}
 is the Serre
functor in ${\cal D}$, in view of the Serre--Grothendieck duality:
$$
{\rm Ext}^i(F\:,\;G)={\rm Ext}^{n-i}(G\:,\;F\otimes\omega_X)^*
$$
for any couple $F, G$ coherent sheaves on $X$ (\cite{Ser}, \cite{Gr}).

For derived categories the translation functor we consider is
 always the usual shift
of grading.

For a closed point $x\in X$ we denote by $k(x)$ the residue field of this point.

We use the standard notations for iterated action of the translation functor
on an object $P$:
$$
P[i] := T^i P, \qquad i\in{\bf Z},
$$
and for the composition of the functors ${\rm Hom}$ and $T$:
$$
{\H i,{}, P, Q}={\rm Ext}^i ( P, Q) := {\h P, {Q[i]}}.
$$

\th{Definition}\label{po}
An object $P\in {\cal D}$ is called {\sf point object} of codimension $s$, if
$$
\begin{array}{ll}
i)& S_{\cal D}(P)\simeq P[s],\\
ii)& {\rm Hom}^{<0}(P\:,\; P)=0,\\
iii)& {\rm Hom}^{0}(P\:,\; P)=k(P).
\end{array}
$$
with $k(P)$ being a field (which is automatically a finite extension of the basic field $k$).
\eth
\th{Proposition}\label{rp}
Let $X$ be a smooth algebraic variety of dimension $n$ with the ample canonical or anticanonical sheaf. Then an object $P\in \db{X}$ is a point object, iff $P\cong{\o x}[r]$, $r\in{\bf Z}$, is isomorphic (up to translation) to the skyscraper sheaf of a closed point $x\in X$.
\eth
\noindent{\sc Remark.} Since $X$ has an ample invertible sheaf it
is projective.
\pr
Any  skyscraper sheaf of a closed point obviously satisfies properties of a point object of the same codimension as the dimension of the variety.

Suppose now that for some object $P\in\db{X}$ properties i)--iii) of definition \ref{po} are verified.

Let ${\cal H}^i$ are cohomology sheaves of $P$. It immediately follows from i) that $s=n$ and ${\cal H}^i\otimes \omega_X={\cal H}^i$. Since $\omega_X$ is  either ample or antiample sheaf, we conclude that ${\cal H}^i$ are finite
length sheaves, i.e. their support are isolated points. Sheaves with the support in different points are homologically orthogonal, therefore any  such object decomposes into direct some of those having the support of all cohomology sheaves in a single point. By iii) the object $P$ is indecomposable. Now consider the spectral sequence, which calculates ${\rm Hom}^m( P\:, \; P)$ by ${\rm Ext}^i( {\cal H}^j\:,\; {\cal H}^k)$:
$$
E^{p,q}_2= \bigoplus_{k-j=q}{\rm Ext}^p({\cal H}^j\:,\; {\cal H}^k) \Longrightarrow {\rm Hom}^{p+q}( P\:, \; P).
$$

Let us mention that for any two finite length sheaves having the same single point as their support, there exists a  non--trivial homomorphism from one to the other, which
sends   generators of the first one to the socle of the second.

Considering ${\rm Hom}^m({\cal H}^j \:, \; {\cal H}^k)$ with minimal $k-j$, we observe that this non--trivial space survives at $E_{\infty}$, hence by ii) $k-j=0$. That means that all but one cohomology sheaves are trivial. Moreover, iii) implies that this sheaf is a  skyscraper. This concludes the proof.

Now having the skyscrapers we are able to reconstruct the invertible sheaves.
\th{Definition}\label{inv}
An {\sf object} $L\in{\cal D}$ is called {\sf invertible} if for any point
object $P\in{\cal D}$ there exists $s\in {\bf Z}$ such that
$$
\begin{array}{lll}
i)& {\rm Hom}^s(L\:,\; P)=k(P),&\\
ii)& {\rm Hom}^i(L\:,\; P)=0, & \quad\mbox{for}\;i\ne s.\\
\end{array}
$$
\eth
\th{Proposition}\label{rin}
Let $X$ be a smooth irreducible algebraic variety. Assume that all point
objects have the form ${\o x}[s]$ for some $x\in X, s\in {\bf Z}$.
Then an object $L\in {\cal D}$ is invertible, iff $L\cong {\cal L}[t]$ for
some invertible sheaf ${\cal L}$ on $X$, $t\in {\bf Z}$.
\eth
\pr
For an invertible sheaf ${\cal L}$ we have:
$$
{\rm Hom}({\cal L}\:, \;{\o x})=k(x), \quad{\rm Ext}^i({\cal L}\:, \; {\o x})=0 ,\quad \mbox{ if }\:i\ne 0.
$$
Therefore, if $L={\cal L}[s]$, then it is an invertible object.

Now let ${\cal H}^i$ are the cohomology sheaves for an invertible object
${L}$. Consider the spectral sequence that calculates
${\rm Hom}^.({ L}\:, \;{\o x})$ for a point $x\in X$ by means of
${\rm Hom}^i({\cal H}^j\:, \;{\o x})$:
$$
E^{p,q}_2 = {\rm Hom}^p({\cal H}^q\:, \;{\o x})\Longrightarrow{\rm Ext}^{p-q}({ L}\:, \;{\o x}).
$$

Let ${\cal H}^{q_0}$ be the nontrivial cohomology sheaf with maximal index. Then for any closed point $x\in X$ from the support of ${\cal H}^{q_0}$ ${\rm Hom}({\cal H}^{q_0}\:, \;{\o x})\ne 0$. But both  ${\rm Hom}({\cal H}^{q_0}\:, \;{\o x})$ and ${\rm Ext}^1({\cal H}^{q_0}\:, \;{\o x})$ are intact by differential of the spectral sequence. Therefore, by definition of an invertible object we obtain that for {\it any} point $x$ from the support of ${\cal H}^{q_0}$
$$
\begin{array}{ll}
a)& {\rm Hom}({\cal H}^{q_0}\:,\;{\o x})=k(x),\\
b)& {\rm Ext}^1({\cal H}^{q_0}\:,\;{\o x})=0.
\end{array}
$$

Since $X$ is smooth and irreducible it follows from b) that  the ${\cal H}^{q_0}$ is 
locally free on $X$. From a) one deduces that  it is invertible.

This implies that ${\rm Ext}^i({\cal H}^{q_0}\:,\;{\o x})=0$  for $i>0$ and
${\rm Hom}({\cal H}^{q_0-1}\:,\;{\o x})$ are intact by differentials of the
spectral sequence. This means that ${\rm Hom}({\cal H}^{q_0-1}\:,\;{\o x})=0$,
 for any $x\in X$, i.e. ${\cal H}^{q_0-1}=0$. Repeating this argument for
${\cal H}^q$ with smaller $q$, we easily see that all ${\cal H}^q$,
except $q=q_0$, are trivial. This proves the proposition.

Now we are ready to prove the reconstruction theorem. Invertible sheaves help us
to `glue' points together.
\th{Theorem}\label{rec}
Let $X$ be a smooth irreducible projective variety with ample canonical or anticanonical sheaf. If ${\cal D}=\db{X}$ is equivalent
  as a graded category to $\db{X'}$ for some other smooth algebraic variety
$X'$, then $X$ is isomorphic to $X'$.
\eth
This theorem is stronger than just a reconstruction for a variety with ample
canonical or anticanonical sheaf from its derived category. Let us mention that since $X'$ might not have ample canonical or anticanonical
sheaf, the situation is not symmetric with respect to $X$ and $X'$.

We divide the proof in several steps, so that the reconstruction procedure
will be transparent.
\pr
During the proof while saying that two isomorphism classes of objects, one in $\db{X}$ and the other in $\db{X'}$, are equal we mean that the former is taken to the latter by the primary equivalence $\db{X}\stackrel{\sim}{\lto} \db{X'}$.

Step 1. Denote ${\cal P}_{ D}$ the set of isomorphism classes of the point objects in ${\cal D}$, ${\cal P}_X$ the set of isomorphism  classes of objects in $\db{X}$
$$
{\cal P}_X:=\Bigl\{ {\o x}[k]\;\Bigl |\: x\in X, k\in {\bf Z}\Bigl\}.
$$
By proposition \ref{rp} ${\cal P}_D\cong{\cal P}_X$. Obviously, ${\cal P}_{X'}\subset {\cal P}_D$. Suppose that there is an object $P\subset {\cal P}_D$, which is not contained in ${\cal P}_{X'}$. Since ${\cal P}_D\cong{\cal P}_X$, any two objects in ${\cal P}_D$ either are homologically mutually orthogonal or belong to a common orbit with respect to the translation functor. It follows that $P\in\db{X'}$ is orthogonal to any skyscraper sheaf
 ${\cal O}_{x'}, x'\in X'$. Hence $P$ is zero. Therefore, ${\cal P}_{X'}=
{\cal P}_D= {\cal P}_X$.

Step 2.  Denote by ${\cal L}_D$ the set of isomorphism classes of invertible
objects in $D$, ${\cal L}_X$ the set of
isomorphism classes of objects in $\db{X}$
$$
{\cal L}_X:=\Bigl\{ L[k]\;\Bigl |\: L\:\mbox{being an invertible sheaf on } X , k\in {\bf Z}\Bigl\}.
$$
By step 1 both varieties $X$ and $X'$ satisfy the assumptions of proposition \ref{rin}. It follows that ${\cal L}_X= {\cal L}_D= {\cal L}_{X'}$.

Step 3.  Let us fix some invertible object $L_0$ in ${\cal D}$, which is an
invertible sheaf on $X$. By step 2 $L_0$ can be regarded, up to translation, as an
invertible sheaf on $X'$. Moreover, changing, if necessary, the equivalence
$\db{X}\simeq\db{X'}$, by the translation functor, we can assume that $L_0$,
regarded as an object on $X'$, is a genuine invertible sheaf.( Formally speaking
$L_0$ is taken by the equivalence $\db{X}\stackrel{\sim}{\lto} \db{X'}$ to
an object which is isomorphic to an invertible sheaf on $X'$. But as was explained
in chapter 1 (formula (\ref{cor})) we can adjust this equivalence so that it takes $L_0$ into the invertible sheaf on $X'$.)

Obviously, by step 1 the set $p_D\subset {\cal P}_D$
$$
p_D:=\Bigl\{ P\in {\cal P}_D \;\Bigl |\: {\h L, P}=k(P)\Bigl\}
$$
coincides with both sets $p_X=\{{\o x}, x\in X\}$ and $p_{X'}=\{{\cal O}_{x'}, x'\in X'\}$. This gives us a pointwise identification of $X$ and $X'$.

Step 4.  Let now $l_X$ (resp., $l_{X'}$) be the subset in ${\cal L}_D$ of isomorphism classes of invertible sheaves on $X$ (resp., on $X'$).

They can be recognized from the graded category structure in $D$ as follows:
$$
l_{X'}=l_X=l_D:=\Bigl\{ L\in {\cal L}_D\;\Bigl |\:{\h L, P}=k(P) \:\mbox{for any } P\in p_D \Bigl\}.
$$
For $\alpha\in{\h L_1, {L_2}}$, where $L_1, L_2\in l_D$, and $P\in p_D$, denote by $\alpha^*_P$ the induced morphism:
$$
\alpha^*_P: {\h L_2, P}\lto {\h L_1, P},
$$
and by $U_{\alpha}$ the subset of those objects $P$ in $p_D$ for which
$\alpha^*_P\ne 0$. By \cite{Il} any algebraic variety has an ample system
of invertible sheaves. This means that $U_{\alpha}$, where $\alpha$ runs over all
elements in ${\h L_1, {L_2}}$ and $L_1$ and $L_2$ runs over all elements in $l_D$,
 give a base for the Zariski topologies on both $X$ and $X'$. It follows that the topologies on $X$ and $X'$ coincide.

Step 5.  Since codimensions of all point objects are equal to the dimensions of $X$ and of $X'$, we have $dimX=dimX'$.
Then, formula (\ref{Ser}) for the Serre functor shows that the operations of twisting by the canonical sheaf on $X$ and on $X'$ induce equal transformations on the set $l_D$.

Let $L_i=F^iL_0[-ni]$. Then $\{L_i\}$ is the orbit of $L_0$ with respect to twisting by the canonical sheaf on $X$. Changing, if necessary, the equivalence $\db{X}\stackrel{\sim}{\lto}\db{X'}$ we can assume that $ \{L_i\}$ is the orbit of $L_0$ with respect to twisting by the canonical sheaf on $X'$ too.

Since the canonical sheaf $\omega_X$ is either ample or antiample, the set of all $U_{\alpha}$, where ${\alpha}$ runs over all elements in ${\h L_i, {L_j}}, i,j\in {\bf Z}$, is the base of the Zariski topology on $X$, hence, by step 4, on $X'$. This means that canonical sheaf on $X'$ is also ample or, respectively, antiample (see \cite{Il}).

For all pairs $(i, j)$ there are natural isomorphisms:
$$
\begin{array}{c}
{\h L_i, {L_j}}\cong{\h S^i L_{0}[-ni], {S^j L_0 [-nj]}}\cong\\
\cong{\h L_0, {S^{j-i}L_0 [-n(j-i)]}}\cong{\h L_0, { L_{j-i}}}.
\end{array}
$$
They induce a ring structure in the graded algebra $A$ over $k$ with graded components
$$
A_i={\h L_0, {L_i}}.
$$

 This algebra , being defined intrinsicly by the graded category structure,
is isomorphic to the coordinate algebra $B$ of the canonical sheaf for $X$,
i.e. algebra with graded components:
$$
B_i={\rm Hom}_{X}({\cal O}_X, \omega^{\otimes i}_{X}).
$$

Indeed, $L_i = L_0\otimes \omega^{\otimes i}_X$, the isomorphism being given by tensoring by $L_0$. It is a ring homomorphism, because the functor of tensoring by $L_0$ commutes with the Serre functor by proposition \ref{com}.

The same is true for the coordinate algebra $B'$ of the canonical sheaf for $X'$.
Eventually, we obtain isomorphism $B\stackrel{\sim}{\to}B'$ of the canonical
algebras on $X$ and $X'$. Since the canonical sheaves on both $X$ and $X'$ are
 ample ( or antiample), both varieties can be obtained by projectivization
from the canonical algebras:
$$
X\stackrel{\sim}{\to}{\bf Proj}B\stackrel{\sim}{\to}{\bf Proj}B'\stackrel{\sim}{\to}X'
$$
This gives a biregular isomorphism between $X$ and $X'$ as algebraic varieties.
 This finishes the proof.


\sec{Group of exact autoequivalences.}

It was explained in chapter 1 that the set of isomorphism classes of exact
autoequivalences in a triangulated category ${\cal D}$ is a group.
We denote this group by ${\rm Aut}{\cal D}$.

The problem of reconstructing of a variety from its derived category is
closely related to the problem of computing the
group of exact autoequivalences for $\db{X}$. For ample canonical or anticanonical sheaf we have the following
\th{Theorem}\label{aut}
Let $X$ be a smooth irreducible projective variety with
ample canonical or anticanonical sheaf. Then the group
of isomorphism classes of exact autoequivalences $\db{X} \to
 \db{X} $ is generated by the automorphisms of the variety, the twists by invertible sheaves and the translations.
\eth

\pr
Assume for definiteness that the canonical sheaf is ample.
Choose an autoequivalence $F$. Since the class of invertible objects is
defined intrinsicly with respect to the graded structure of the category,
it is preserved by any autoequivalence. Moreover, the set of isomorphism
classes of invertible objects is transitive with respect to the action of
the subgroup ${\rm Aut}\db{X}$ generated by translations and twists. Indeed, by propositions \ref{rp} and \ref{rin} all invertible objects in $\db{X}$ are invertible sheaves up to translations. Any invertible sheaf can be obtained from the trivial sheaf ${\cal O}$ by applying the functor of tensoring with this invertible sheaf. Therefore, using twists with invertible sheaves and translations we can assume that our functor $F$ takes ${\cal O}$ to ${\cal O}$. It follows that $F$ takes any tensor power $\omega^{\otimes i}_X$ of the canonical sheaf into itself, because by proposition \ref{com} it commutes with the Serre functor.

Therefore, our functor induces an automorphism of the graded coordinate algebra $A$ of the canonical sheaf , i.e. algebra with graded components:
$$
A_i={\h {\cal O}, {\omega^{\otimes i}_X}}={\rm H}^0 (\omega^{\otimes i}_X)
$$

Any graded automorphism of the canonical algebra induces an automorphism of the variety. Adjusting our functor $F$ by an autoequivalence induced by an automorphism of the variety we can assume that the automorphism of the canonical algebra induces the trivial automorphism of the variety. 

Such an automorphism is actually a scaling, i.e. it takes an element $a\in {\rm H}^0 (\omega^{\otimes i}_X)$ to $\lambda^i a$, for some fixed scalar $\lambda$. Indeed, the graded ideal generated by any element $a\in {\rm H}^0 (\omega^{\otimes i}_X)$ is stable with respect to the automorphism. It follows that $a$ is multiplied by a scalar. Then the linear operator in the graded component ${\rm H}^0 (\omega^{\otimes i}_X)$ induced by the automorphism  should be scalar , say $\lambda_i$. Since our automorphism is algebra automorphism, it follows that $\lambda_i =\lambda^i$, for $\lambda = \lambda_1$ (in case ${\rm H}^0 (\omega_X )=0$, i.e. when $\lambda_1$ is not
defined, we may substitute in this reasoning the Serre functor by a sufficient  $j$-th power of it such that ${\rm H}^0 (\omega_X^{\otimes j})\ne 0$, and respectively the canonical algebra by the corresponding Veronese subalgebra).

To kill the scaling of the canonical algebra we substitute functor $F$ by an isomorphic one. For this we take the subclass ${\cal C}$ of objects in ${\cal D}$ consisting of powers of the canonical sheaf:
$$
{\cal C}=\{ \omega^{\otimes i}_X \},\qquad i\in{\bf Z}.
$$
As in chapter 1 for any object $C$ in ${\cal C}$ we need to choose an
isomorphism of its image with some other object $C'$. Our functor preserves
all objects from ${\cal C}$. We choose $C'=C$ for any $C$ in ${\cal C}$ and
the non--trivial isomorphism: if $C=\omega^{\otimes i}_X$, then the isomorphism is $\lambda^{-i}\cdot id_C$.

Then the new functor constructed by formula (\ref{cor}) induces the trivial
automorphism of the canonical ring.

Thus we have a functor, which takes the trivial invertible sheaf and any power of the canonical sheaf to themselves and preserves homomorphisms between all these sheaves. Let us show that such a functor is isomorphic to the identity functor.

First, our functor takes pure sheaves to objects, isomorphic to pure sheaves,
because such objects can be characterized as the objects ${\cal G}$ in
$\db{X}$ having trivial ${\rm Hom}^k (\omega^{\otimes i}_X , {\cal G})$,
for $k\ne 0$ and for sufficiently negative $i$.
Again we can substitute our functor by an isomorphic one, which takes sheaves to pure sheaves.
By Serre theorem \cite{Ser} the abelian category of pure sheaves is equivalent to the category of graded finitely generated modules over the canonical algebra $A$ modulo the subcategory of finite dimensional modules. The equivalence takes a sheaf ${\cal G}$ into a module ${\cal M}({\cal G})$ with graded components:
$$
{\cal M}_i ({\cal G})={\h \omega^{\otimes -i}_X , {\cal G}}
$$
Our functor $F$ gives isomorphisms:
$$
{\h \omega^{\otimes -i}_X , {\cal G}}\stackrel{\sim}{\lto}{\h F(\omega^{\otimes -i}_X) , {F({\cal G})}}\cong{\h \omega^{\otimes -i}_X , {F({\cal G})}}
$$

Since $F$ induces trivial action on the canonical algebra these isomorphisms form  an isomorphism of $A$-modules:
$$
{\cal M}({\cal G})\stackrel{\sim}{\to}{\cal M}(F({\cal G})).
$$

It is natural with respect to ${\cal G}$. Hence we obtain an isomorphism of functors ${\cal M}\stackrel{\sim}{\to}{\cal M}\circ F$.

Since modulo the subcategory of finite dimensional modules  ${\cal M}$ is
an equivalence, we have a functorial isomorphism $id\stackrel{\sim}{\to}F$
on the subcategory of coherent sheaves.

Our system of objects $\{\omega^{\otimes i}_X \}$ enjoys nice properties
with respect to the abelian category of coherent sheaves on $X$, which
allow to extend the natural transformation $id\to F$, from the core
of the t-structure to a natural isomorphism in the whole derived category. It was done in Proposition \ref{ext} of the Appendix.

This finishes the proof of the theorem.

 In the hypothesis of Theorem \ref{aut} the group ${\rm Aut}{\db{X}}$ is the semi-direct product of its subgroups $G_1 = {\rm Pic} X \oplus {\bf Z}$ and $G_2 = {\rm Aut} X$, ${\bf Z}$ being generated by the translation functor:
$$
{\rm Aut}{\db{X}}\cong {\rm Aut} X \ltimes ({\rm Pic} X \oplus {\bf Z})
$$
 Indeed, in course of the proof of
the theorem we in fact showed that any element from ${\rm Aut}{\db{X}}$ could be decomposed as $g=g_1 g_2$ with $g_1 \in G_1$ and $g_2 \in G_2$.
The subgroups $G_1$ and $G_2$ meet trivially in $G$, because the elements from the latter take the structure sheaf ${\cal O}$ to itself, while those from the former do not.
Group $G_1$ is obviously preserved by conjugation by elements from $G_1$
and $G_2$, hence normal in $G$.

\sect{Appendix.}

This appendix is devoted to describing conditions, under which one can extend to the whole category a natural isomorphism between the identity functor and an exact autoequivalence in the bounded derived category ${\db{\cal A}}$, provided one has such an isomorphism in an abelian category ${\cal A}$ (or even in a smaller subcategory, see the proposition below).

To break our way through technical details we need a sequence of objects in the abelian category with some remarkable properties. For the case when the sequence consists of powers of an invertible sheaf these properties are resulted from ampleness of this sheaf. For this reason we postulate them under the name of ampleness.

\ap{Definition}\label{obil}
Let ${\cal A}$ be an abelian category. We call a sequence of objects
$\{ P_i \}, \quad i\in {\bf Z}_{\le 0}$,  {\sf ample} if for every object $X\in {\cal A}$,
there exists $N$  such that for all $i< N$ the following conditions hold:

a) the canonical morphism ${\h P_i, X}\otimes P_i \lto X$ is surjective,

b) ${\E j, {}, P_i, X}=0$ for any $j\ne 0$,

c) ${\h X, {P_i}}=0$.
\eap

Denote by $\d{\cal A}$  the bounded derived category of ${\cal A}$. Let us
consider ${\cal A}$ as a full subcategory $j: {\cal A}\hookrightarrow\d{\cal A}$
in $\d{\cal A}$ in the usual way. We also consider a full subcategory $q : {\cal C}\hookrightarrow {\d{\cal A}}$ with $Ob{\cal C}=\{ P_i \}_{i\in {\bf Z}_{\le 0}}$. We
shall show that if there exists an exact autoequivalence
$F : \d{{\cal A}}\lto\d{{\cal A}}$ and an isomorphism of its restriction to
${\cal C}$ with the identity functor $id_{\cal C}$, then this isomorphism
 can be uniquely extended to an isomorphism of $F$ with the identity functor
$id_{\d{\cal A}}$ in  the whole $\d{\cal A}$. 

The idea is in reducing
the number of cohomology for an object by killing the highest one by
means of a surjective morphism from $\oplus P_i$ for sufficiently negative $i$.

For the proof we shall repeatedly use the following lemma (see \cite{BBD}).

\ap{Lemma}\label{tr}
Let $g$ be a morphism from $Y$ to $Y'$ and suppose that these objects
  are included into
the following two exact triangles:
$$
\begin{array}{ccccccc}
X&\stackrel{u}{\lto}&Y&\stackrel{v}{\lto}&Z&\stackrel{w}{\lto}&X[1]\\
\dda\rlap{\ss{f}}&&\da\rlap{\ss{g}}&&\dda\rlap{\ss{h}}&&\dda\rlap{\ss{f[1]}}\\
X'&\stackrel{u'}{\lto}&Y'&\stackrel{v'}{\lto}&Z'&\stackrel{w'}{\lto}&X'[1]
\end{array}
$$
If $v'gu=0$, then there exist morphisms $f : X\to X'$ and $h : Z\to Z'$ such
that the triple $(f, g, h)$ is a morphism of triangles.

If, in addition, ${\h X[1], {Z'}}=0$, then the morphisms $f$ and $h$,
making commutative the first   and, respectively, the second square of the diagram, are
 unique.
\eap

\ap{Proposition}\label{ext}
Let ${\cal A}$ be an abelian category possessing an ample sequence
$\{P_i \}$ and let $F :\d{{\cal A}}\lto\d{{\cal A}}$ be an exact
autoequivalence.
Suppose there exists an isomorphism
$f : q\stackrel{\sim}{\lto}F\mid_{\cal C} $ (where $q : {\cal C}\hookrightarrow
\d{\cal A}$ is the natural embedding).
 Then this isomorphism can be uniquely extended to an isomorphism
$id\stackrel{\sim}{\lto}F$ in the whole $\d{\cal A}$.
\eap
\pr
Note that $X\in {\d{\cal A}}$ is isomorphic to an object in ${\cal A}$ iff ${\rm Hom}^j (P_i , X )=0$ for $j\ne 0$ and $i\ll 0$, in view of the condition b) from Definition \ref{obil}.

This allows us to "extract" the abelian subcategory  ${\cal A}$ from ${\d{\cal A}}$ by means of the sequence $\{P_i \}$. 
Then using surjective coverings ${\h P_i, X}\otimes P_i \lto X$ by standard techniques from the theory of abelian categories one can extend $f$ to an isomorphism (which we denote by the same letter) $f : j\stackrel{\sim}{\lto} F \mid_{\cal A}$, where $j$ stands for the natural embedding $j: {\cal A}\hookrightarrow\d{\cal A}$. We skip details of this part of the proof, because
we don't need it in the main body of the paper.

Let us define $f_{X[n]} : X[n]\lto F(X[n])\cong F(X)[n]$ for $X\in{\cal A}$
by
$$
f_{X[n]}=f_X [n].
$$
It is not difficult to show that for any $X$ and $Y$ in ${\cal A}$ and for any
$u\in{\E k, {}, X, Y}$ the diagram
\begin{equation}\label{duf}
\begin{array}{ccc}
X&\stackrel{u}{\lto}&Y[k]\\
\llap{\ss{f_X}}\da&&\da\rlap{\ss{f_Y [k]}}\\
F(X)&\stackrel{F(u)}{\lto}&F(Y)[k].\\
\end{array}
\end{equation}
is commutative.
Indeed, since any element $u\in{\E k, {}, X, Y}$ can be represented
as the Yoneda composition $u = u_1\cdots u_k$
of elements  $u_i\in{\E 1, {}, Z_i, {Z_{i+1}}}$ for some
objects $Z_i$, with $Z_1=X, Z_{k+1}=Y$,
we can restrict ourselves to the case $u\in{\E 1, {}, X, Y}$.
Consider  the following diagram:
$$
\begin{array}{ccccccc}
Y&\lto&Z&\stackrel{p}{\lto}&X&\stackrel{u}{\lto}&Y[1]\\
\llap{\ss{f_Y}}\da&&\da\rlap{\ss{f_Z}}&&\dda\rlap{\ss{h}}&&\da\rlap{\ss{f_Y [1]}}\\
F(Y)&\lto&F(Z)&\stackrel{F(p)}{\lto}&F(X)&\stackrel{F(u)}{\lto}&F(Y)[1]\\
\end{array}
$$
By an axiom of triangulated categories there exists a morphism
$h : X\to F(X)$
such that $(f_Y, f_Z, h)$ is a morphism of triangles.
On the other hand, since ${\h Y[1], {F(X)}}=0$, by the lemma above $h$
is a unique morphism
such that  $F(p)\circ f_Z = h\circ p$.
As $F(p) f_Z = f_X p$, we conclude that $h=f_X$. This implies
commutativity of the diagram (\ref{duf}) for $k=1$.

We shall prove by induction over $n$ the following statement.
Consider the full subcategory $j_n : {\cal D}_n \hookrightarrow \d{\cal A}$
in $\d{\cal A}$ generated by objects having nontrivial cohomology in a (non-fixed) segment
of length $n$. Then there is a unique extention of $f$ to a natural
functorial isomorphism $f_n : j_n \lto F\mid_{{\cal D}_n}$.

Above we have completed the first, $n=1$, step of the induction.

Now take the step $n=a, \; a\ge 1$,  for granted.
Let $X$ be  an object in $ {\cal D}_{a+1}$ and suppose, for definiteness,
 that its cohomology
${\cal H}^p (X)$ are nontrivial only for $p\in [-a, 0]$. Take
$P_i$ from the given ample sequence with sufficiently negative $i$ such that
\begin{equation}
\begin{array}{ll}\label{tt}
a)& {\H j, {}, P_i, {{\cal H}^p (X)}}=0 \mbox{ for all }\; p\; \mbox{ and for }
\; j\not=0,\\
b)&\mbox{there exists a surjective morphism}\quad u: P_i^{\oplus k}\lto {\cal H}^0 (X),\\
c)&{\h {\cal H}^0 (X), {P_i}}=0.
\end{array}
\end{equation}
Note that in view of condition a) and the standard spectral sequence ${\h P_i , X}\stackrel{\sim}{\to}
{\h P_i , {{\cal H}^0 (X)}}$. This means that we can find a morphism
$v: P^{\oplus k}_i \lto X$ such that the composition of
$v$ with the canonical morphism $X\lto {\cal H}^0 (X)$ coincides with $u$.
Consider an exact triangle:
$$
Y[-1]\lto P^{\oplus k}_i\stackrel{v}{\lto}X\lto Y.
$$

Denote by $f_i$ the morphism $f_X$ for $X=P^{\oplus k}_i$.
Since $Y$ belongs to ${\cal D}_a$ by the induction hypothesis,  the isomorphism
$f_Y$ is already defined  and  the diagram:
\begin{equation}
\begin{array}{ccccccc}\label{dnt}
P^{\oplus k}_i&\stackrel{v}{\lto}&X&\lto&Y&\lto&P^{\oplus k}_i [1]\\
\da\rlap{\ss{f_i}}&&\dda\rlap{\ss{f_X}}&&\da\rlap{\ss{f_Y }}&&\da\rlap{\ss{f_i [1]}}\\
F(P^{\oplus k}_i)&\stackrel{F(v)}{\lto}&F(X)&\lto&F(Y)&\lto&F(P^{\oplus k}_i)[1]\\
\end{array}
\end{equation}
is commutative.

Further,  we have the following sequence of isomorphisms:
$$
{\h X, {F(P^{\oplus k}_i)}}\cong{\h X, {P^{\oplus k}_i}}\cong
{\h {\cal H}^0 (X), {P^{\oplus k}_i}}=0.
$$
Hence, applying  lemma \ref{tr} to $g$ equal $f_Y$, we obtain a unique
morphism $f_X : X\lto F(X)$ that preserves commutativity of the above diagram.

It is clear from the definition that $f_X$ is an isomorphism, if so are
$f_i$ and $f_Y$.
For the future we need to show that $f_X$ does not depend on the choice for
$i$ and $u$.
Suppose we are given two surjective morphisms $u_1 : P_{i_1}^{\oplus k_1}\lto {\cal H}^0 (X)$
and $u_2 : P_{i_2}^{\oplus k_2}\lto {\cal H}^0 (X)$, where $i_1$ and
$i_2$ are sufficiently negative to satisfy conditions a), b) and c).
Then we can find sufficiently negative $j$ and
surjective morphisms $w_1, w_2$ such that the following diagram  commutes:
$$
\begin{array}{ccc}
P^{\oplus l}_{j}&\stackrel{w_2}{\lto}& P^{\oplus k_2}_{i_2}\\
 \da\rlap{\ss{w_1}}&& \da\rlap{\ss{u_2}}\\
P^{\oplus k_1}_{i_1}&\stackrel{u_1}{\lto}&{\cal H}^0 (X).\\
\end{array}
$$

Denote by $v_1 : P_{i_1}^{\oplus k_1}\lto  X, v_2 : P_{i_2}^{\oplus k_2}\lto  X$
the morphisms corresponding to $u_1$ and $u_2$.
Since ${\h P_j , X}\stackrel{\sim}{\to}{\h P_j , {{\cal H}^0 (X)}}$,
we have $v_2 w_2 = v_1 w_1$.

There is a morphism $\phi : Y_j\lto Y_{i_1}$ such that  the triple
$(w_1, id, \phi)$ is a morphism of exact triangles:
$$
\begin{array}{ccccccc}
P^{\oplus l}_j&\stackrel{v_1\circ w_1}{\lto}&X&\stackrel{y}{\lto}&Y_j&\lto&P^{\oplus l}_j [1]\\
\llap{\ss{w_1}}\da&&\da\rlap{\ss{id}}&&\da\rlap{\ss{\phi}}&&\da\rlap{\ss{w_1 [1]}}\\
P^{\oplus k_1}_{i_1}&\stackrel{v_1}{\lto}&X&\stackrel{y_1}{\lto}&Y_{i_1}&\lto&
 P^{\oplus k_1}_{i_1} [1],\\
\end{array}
$$
i.e. $\phi y = y_1$.

Since $Y_j$ and $Y_{i_1}$ have cohomology in the segment $[-a, -1]$,
by the induction hypothesis, the following square is commutative:
$$
\begin{array}{ccc}
Y_j&\stackrel{\phi}{\lto}&Y_{i_1}\\
\llap{\ss{f_{Y_j}}}\da&&\da\rlap{\ss{f_{Y_{i_1}}}}\\
F(Y_j)&\stackrel{F(\phi)}{\lto}&F(Y_{i_1}).\\
\end{array}
$$

Denote by $f^j_X , f^{i_1}_X , f^{i_2}_X$ the unique morphisms constructed as
above to make commutative the diagram (\ref{dnt}) for $v$ equal
respectively $v=v_1 w_1 , \; v=v_1 , \; v=v_2$.
Further, we have:
$$
F(y_1 ) f^j_X = F(\phi y) f^j_X = F(\phi ) F(y) f^j_X = F(\phi ) F_{Y_j} y =
f_{Y_{i_1}} \phi y = f_{Y_{i_1}} y_1 .
$$
It follows that $f^j_X = f^{i_1}_X$. Analogously, since $v_1 w_1 = v_2 w_2$
we have $f^j_X = f^{i_2}_X$. Therefore the morphism $f_X$ does not
depend on the choice of $i$ and of the morphism $u: P^{\oplus k}_i\lto {\cal H}^0 (X)$.

By means of the translation functor we obtain in the obvious way
the only possible extension of $f_a$ to ${\cal D}_{a+1}$. Let us prove
that it is indeed a natural transformation from $j_{a+1}$ to $F\mid_{{\cal D}_{a+1}}$
, i.e. that for any morphism $\phi : X\lto Y$, $X, Y$ being in ${\cal D}_{a+1}$,
 the following
 diagram commutes

\begin{equation}\label{comd}
\begin{array}{ccc}
X&\stackrel{\phi}{\lto}&Y\\
\llap{\ss{f_X}}\da&&\da\rlap{\ss{f_Y }}\\
F(X)&\stackrel{F(\phi)}{\lto}&F(Y).
\end{array}
\end{equation}

We shall reduce the problem to the case when both  $X$ and $Y$ are
in ${\cal D}_a$.

 There are two working posibilities that we shall utilize for this.

{\it Case 1.}
Suppose that the upper bound, say $0$ (without loss of generality), of
cohomology for $X$ is greater than that for $Y$.
Take a surjective morphism $u : P^{\oplus k}_i\lto{\cal H}^0 (X)$
satisfying a), b), c) and construct the morphism
 $v: P^{\oplus k}_i \lto X$ related to $u$ as above.
We have an exact triangle:
$$
\begin{array}{ccccccc}
P^{\oplus k}_i&\stackrel{v_1}{\lto}&X&\stackrel{\alpha}{\lto}&Z&\lto&
P^{\oplus k}_i [1].\\
\end{array}
$$
If we take $i$ sufficiently negative, then ${\h P^{\oplus k}_i, Y}=0$.
Applying the functor ${\h {-}, Y}$ to this triangle we found that
there exists a morphism $\psi : Z\lto Y$ such that $\phi=\psi \alpha$.
We know that $f_X$, defined above, satisfies the equation
$$
F(\alpha ) f_X = f_Z \alpha .
$$

If we assume that
$$
F(\psi ) f_Z = f_Y \psi,
$$
then
$$
F(\phi ) f_X = F(\psi ) F(\alpha ) f_X = F(\psi ) f_Z \alpha = f_Y \psi \alpha =
f_Y \phi .
$$

This means that for this  case in verifying  commutativity of (\ref{comd})
we can substitute $X$ by an object $Z$ such that the upper bound of
its cohomology is less by one than that for $X$. Moreover, one can easily see that
if $X$  belongs to ${\cal D}_k$, with $k>1$, then $Z$ does to ${\cal D}_{k-1}$,
and if it is in ${\cal D}_1$ then so is $Z$.

{\it Case 2.}
Suppose now that the upper bound , say $0$ (again without loss of generality),
of cohomology for $Y$ is greater than or equal to that for $X$.
Take a surjective morphism $u : P^{\oplus k}_i\lto{\cal H}^0 (Y)$ with
$i$ satisfying a), b), c) (with $Y$ instead of $X$) and construct a morphism
$v: P^{\oplus k}_i\lto Y$ related to $u$.
Consider an exact triangle
$$
\begin{array}{ccccccc}
P^{\oplus k}_i &\stackrel{v}{\lto}&Y&\stackrel{\beta}{\lto}&W&\lto&
P^{\oplus k}_i.\\
\end{array}
$$

By $\psi$ denote the composition $\beta\circ\phi$.

If we assume that
$$
F(\psi ) f_X = f_W \psi,
$$
then, since $F(\beta ) f_Y = f_W \beta$ we have:
\begin{equation}\label{cf}
F(\beta)(f_Y \phi - F(\phi) f_X) = f_W \beta \phi - f(\beta\phi ) f_X =
f_W \psi - F(\psi ) f_X=0.
\end{equation}

We again take $i$ sufficiently negative, so that  ${\h X, {P^{\oplus k}_i }}=0$.
As $F(P^{\oplus k}_i)$ is isomorphic to $P^{\oplus k}_i$, then
${\h X, {F(P^{\oplus k}_i )}}=0$.
Applying the functor ${\h X, {F(-)}}$ to the above triangle we found that the
composition with $F(\beta )$ gives an inclusion of ${\h X, {F(Y)}}$ into
${\h X, {F(W)}}$.
It follows from (\ref{cf}) that
$f_Y \phi = F(\phi) f_X$.

Thus in this case in verifying  commutativity of (\ref{comd})
we can substitute $Y$ by an object $W$ such that the upper bound of
its cohomology is less by one than that for $Y$.
If $Y$ belongs to ${\cal D}_k , \; k>1$, then $W$ does to ${\cal D}_{k-1}$,
if $Y$ belongs to ${\cal D}_1$, then so does $W$.

Suppose now that $X$ and $Y$ are in ${\cal D}_{a+1} , \; a>1$. Depending
on  which case, 1) or 2), we are in, we can substitute either $X$ or $Y$ by an
object lying in ${\cal D}_a$. Then repeating, if necessary, the procedure
we can lower the upper bound of the cohomology of the object to such
a point that the other case is applicable. Then we shorten the cohomology segment
of the second object and come to the situation when both objects are in
${\cal D}_a$, i.e. to the induction hypothesis.

At every step of the construction we always made the only possible choice
for the  morphism $f_{X}$. This means that the natural transformation  with
required properties is unique.

This finishes the proof of the proposition.

Algebra section,
Steklov Math. Institute, Vavilova 42, 117333 Moscow, RUSSIA

\noindent bondal@sovam.com , \quad orlov@class.mi.ras.ru
\end{document}